\documentclass[12pt]{iopart}
\usepackage{iopams}
\usepackage{perpage} %the perpage package
\MakePerPage{footnote} %the perpage package command
\usepackage{graphics}
\usepackage{relsize}
\usepackage{latexsym}
\usepackage{color}
\usepackage{epstopdf}
\usepackage{ctable}
\usepackage{graphicx}
\usepackage{array,multirow,makecell}
\usepackage{verbatim}
\newtheorem{theo}{Theorem}[section]
\newtheorem{uppg}[theo]{Exercise}
\newcommand{\be}{\begin{eqnarray*}}
\newcommand{\ee}{\end{eqnarray*}}
\newcommand{\ben}{\begin{eqnarray}}
\newcommand{\een}{\end{eqnarray}}

\def\subsecn (#1) {\medskip\ \ \ {\it #1}\medskip}

\newcommand{\lp}[1]{\left(\begin{array}{#1}}
\newcommand{\rp}{\end{array}\right)}

\newcommand{\leftd}[1]{\left\{\begin{array}{#1}}
\newcommand{\rightd}{\end{array}\right.}
\setcellgapes{1pt}
\makegapedcells
\newcolumntype{R}[1]{>{\raggedleft\arraybackslash }b{#1}}
\newcolumntype{L}[1]{>{\raggedright\arraybackslash }b{#1}}
\newcolumntype{C}[1]{>{\centering\arraybackslash }b{#1}}

\def\A {\mathbf{A}}

\def\I {\mathbf{I}}

\def\R {\mathbf{R}}

\def\b {\boldsymbol{b}}

\def\e {\boldsymbol{e}}
\def\f {\boldsymbol{f}}

\def\h {\boldsymbol{h}}

\def\k {\boldsymbol{k}}
\def\l {\boldsymbol{l}}
\def\m {\boldsymbol{m}}
\def\n {\boldsymbol{n}}

\def\u {\boldsymbol{u}}
\def\v {\boldsymbol{v}}
\def\w {\boldsymbol{w}}
\def\x {\boldsymbol{x}}
\def\y {\boldsymbol{y}}
\def\z {\boldsymbol{z}}

\def\Nc {\mathcal{N}}

\def\Rb {\mathbb{R}}
\def\Pb {\mathbb{P}}
\def\Eb {\mathbb{E}}

\date{\today}
\begin{document}

\title[]{Multilevel Monte Carlo simulation of a diffusion
with non-smooth drift.}

\author{Azzouz Dermoune$^1$\footnote{Corresponding 
author: Azzouz.Dermoune@univ-lille1.fr}, Daoud Ounaissi$^2$ and Nadji Rahmania$^3$}

\address{$^1$Dermoune Azzouz, Cit\'e scientifique, France}
\address{$^2$Daoud Ounaissi, Cit\'e scientifique, France}
\address{$^3$Nadji Rahmania, Cit\'e scientifique, France}
\eads{\mailto{Azzouz.Dermoune@univ-lille1.fr}, \mailto{daoud.ounaissi@ed.univ-lille1.fr}, \mailto{nadji.rahmania@univ-lille1.fr}}
% For more than one e-mail address, please use the command \eads{\mailto{#1}, \mailto{#2}} with \mailto surrounding each e-mail address.

\begin{abstract}
We show that Lasso and Bayesian Lasso are very close when the sparsity is large
and the noise is small. Then we propose to solve Bayesian Lasso using multivalued stochastic
differential equation. We obtain three discretizations algorithms, and propose a method
for calculating the cost of Monte-Carlo (MC), multilevel Monte Carlo (MLMC) and MCMC
algorithms.
\end{abstract}

\vspace{2pc}

\section{Introduction}
%\section{Motivation} 
Let $\y=\A\x+\sigma\w$ be the classical linear regression problem see e.g. 
\cite{Tibshirani1996} and the references herein, 
(see also \cite{DRW,DDR1,DDR2} for some new applications). 
Here $p$ and $n$ is a couple of positive integers, $\y\in\Rb^n$ are the observations, $\x\in \Rb^p$ is the unknown signal to recover, $\w\in\Rb^n$ is the standard noise, $\sigma$ is the size of the noise and $\A$ is a known matrix which maps the signal domain $\Rb^p$ into the observation domain $\Rb^n$.  
The matrix $\A$ is in general ill-conditioned (e.g. in the case $n < p$)  
which makes difficult to use the least squares estimate. Penalization is a popular way to compute an approximation of $\x$ from 
the observations $\y$. The general framework proposes to recover the vector $\x$ using 
the posterior probability distribution function proportional to 
\be 
\exp\left(-P(\x)-\frac{\|\A\x-\y\|^2}{2\sigma^2}\right).
\ee 
Here $\|\cdot\|$ denotes the Euclidean norm.     
This requires to define a penalization $P$ to enforce some prior information on 
the signal $\x$. The term $\frac{\|\A\x-\y\|^2}{2}$ reflects Gaussian prior on the noise $\w$. The parameter $\sigma^2>0$ reflects the noise level. 

The $l^1$ penalization is the sum of the absolute values $P(\x)=\alpha \|\x\|_1$ of the components of $\alpha\x$. The parameter $\alpha >0$ reflects the sparsity level of the variable $\x$. 
The Lasso := $\arg\min\{\alpha \|\x\|_1 + \frac{\|\A\x-\y\|^2}{2}, \x \in \Rb^p\}$ was first introduced in \cite{Tibshirani1996}. It is also called Basis Pursuit De-Noising method \cite{Donoho}.  
It was introduced to induce sparsity in the variable $\x$. A large number of theoretical results has been provided 
for the $l^1$ penalization see e.g. \cite{Daubechies2004,Dossal,vaiter} and the references 
herein.   

We will suppose that $\alpha=2\beta$ and $\sigma^2=\frac{1}{2\beta}$. 
It follows that the posterior PDF is equal to  
\ben 
\frac{1}{Z_\beta}\exp\left(-2\beta F(\x)\right), 
\label{posterior}
\een 
where 
\ben
F(\x)=\|\x\|_1+\frac{\|\A\x-\y\|^2}{2}
\label{F}
\een
and $Z_\beta$ is the partition function, i.e. 
\be 
Z_\beta=\int_{\Rb^p}\exp\big(-2\beta F(\x)\big)d\x.
\ee 

Bayes estimator of $\x$ is equal to 
\ben 
\m_\beta:=\int_{\Rb^p}\x\exp\big(-2\beta F(\x)\big)\frac{d\x}{Z_\beta}.
\label{bayes} 
\een 
Lasso is the maximum a posteriori estimator
\ben 
Lasso =\arg\min\big\{F(\x):\quad \x\in\Rb^p\big\}.
\label{Lasso} 
\een 
In the sequel $X_\beta$ will denote a random vector having the probability distribution (\ref{posterior}). Hence Bayes estimator (\ref{bayes}) 
is the mathematical expectation 
\ben 
\Eb[ X_\beta ]. 
\label{mathematicalexpectation} 
\een  

In the first part of this work we show how Bayes estimator converges to Lasso as $\beta\to +\infty$.
In the second part we consider for fixed $\beta$ the random vector $X_\beta$ 
as the limit of a multivalued stochastic process $(\x(T))$ (Langevin diffusion with non-smooth drift) 
as $T\to +\infty$. 
We propose to approximate Bayes estimator $\Eb[ X_\beta ]$ by 
the mathematical expectation $\Eb[\x(T)]$ for large $T$.  
We obtain three discretizations algorithms. 
Two among them are known as
unadjusted Langevin algorithm (ULA) (\cite{Pereyra}) and STMALA (\cite{Fort}). 
We calculate the latter mathematical expectation $\Eb[\x(T)]$ using Monte Carlo (MC),  
Multilevel Monte Carlo (MLMC) and MCMC methods.  
We propose a method for calculating the cost of MC, MLMC and MCMC.  
\section{Lasso estimator properties}  
First, we need some notations. For each 
$\x\in\Rb^p$, the sub-differential 
$sgn(\x)=\partial\|\x\|_1$ is the set of  
the column vector $\xi\in\Rb^p$ such that  the component 
$\xi_i=sgn(x_i)=1$ if $x_i >0$, $\xi_i=sgn(x_i)=-1$ if $x_i <0$ and $\xi_i\in [-1,1]$ if $x_i=0$.  
 
We will denote, for each subset $J\subset\{1, \ldots, p\}$ and for each vector $\v\in \Rb^p$, $\v(J)=\left(v(i): i\in J\right)\in \Rb^{|J|}$. Here $|J|$ denotes the cardinality of $J$. 
The notation $\v\leq \w$ means $v(i)\leq w(i)$ for all $i=1, 2, \ldots, p$. The scalar product is denoted by $\langle\cdot, \cdot\rangle$, 
and $(\e_i:\quad i=1, 2, \ldots, p)$ denotes the canonical basis of $\Rb^p$.

Now we recall a well known properties of Lasso estimator see e.g. \cite{Tibshirani2013}. 

{\bf lemma:} \label{Lassoproperties} The vector $\x(\y)$ is a minimizer of  the map 
$\x\to F(\x)=\|\x\|_1+\frac{\|\A\x-\y\|^2}{2}$ if and only if 
the vector 
\ben 
\xi:=\A^*(\y-A\x(\y))\in sgn(\x(\y)). 
\label{xi} 
\een 
The vectors  
$\xi$, $\A\x(\y)$ and the $l^1$-norm $\|\x(\y)\|_1$ are constant on the set of Lasso estimators. Moreover, the set of Lasso is convex and compact.  
Here $\A^*$ denotes the transpose of the matrix $\A$. \\

We introduce the sets 
\ben 
I&=&\big\{i\in\{1,\ldots, p\}:\quad |\xi_i| < 1\big\},\label{I}\\ 
\partial I&=&\big\{i\in \{1,\ldots, p\}:\quad |\xi_i|=1\big\}. \label{partialI}
\een 
Observe that the support $\{i\in\{1, \ldots, p\}:\quad x_i(\y)\neq 0\}$ of any Lasso $\x(\y)$ is contained in $\partial I$, and $I$ is contained in the set 
$\{i\in\{1,\ldots, p\}: \quad x_i(\y)=0\}$
of the null components
of $\x(\y)$. For each subset $J$ of $\{1, \ldots, p\}$, 
$\A_J$ denotes the submatrix of $\A$ having its columns indexed by $J$. 

From "equation~(\ref{xi})" it is easy to show that the injectivity of $\A_{\partial I}$  
implies the uniqueness of Lasso. In fact, under this hypothesis the system 
\be 
\xi_{\partial I}=\A_{\partial I}^T\y-\A_{\partial I}^T\A_{\partial I}\x_{\partial I}(\y)
\ee 
has a unique solution. As the support of any Lasso $\x(\y)$ is contained in $\partial I$, then 
Lasso is unique.

In the sequel for each $\x\in\Rb^p$,  
\be 
\pi(\x)=\arg\min\{\|\x-\x(\y)\|:\quad \x(\y)\in Lasso\}.  
\ee 
{\bf prop:}\label{CVP} 
The random positive number $\|X_\beta-\pi(X_\beta)\|$ converges to 0 in 
probability as $\beta\to +\infty$. \\

{\bf proof:} The proof is similar to Theorem 4.1. in  \cite{AthreyaHwang}.
It works as following. 
  
Let $\delta >0$, and $\eta>0$ such that 
\be 
\inf\{F(\x):\quad \|\x-\pi(\x)\|\geq \delta\}> M(\eta)=\sup\big\{F(\x):\quad \|\x-\pi(\x)\|\leq \eta\big\},
\ee  
where $F$ is given by "equation~(\ref{F})". 
We have 
\be 
P(\|X_\beta-\pi(X_\beta)\|\geq \delta)&=&
\frac{\int_{\|\x-\pi(\x)\|\geq \delta}\exp(-\beta F(\x))d\x}
{\int\exp(-\beta F(\x))d\x}\\
&\leq & \frac{\int_{\|\x-\pi(\x)\|\geq \delta}\exp\big(-\beta (F(\x)-M(\eta)\big)d\x}
{\int_{\|\x-\pi(\x)\|\leq \eta}\exp\big(-\beta (F(\x)-M(\eta)\big)d\x}.  
\ee    
From the estimate 
\be 
\int_{\|\x-\pi(\x)\|\geq \delta}\exp\big(-\beta (F(\x)-M(\eta))\big)d\x\leq 
\int_{\|\x-\pi(\x)\|\geq \delta}\exp\big(-(F(\x)-M(\y))\big)dx < +\infty
\ee 
and the bounded convergence theorem, the numerator 
$\int_{\|x-\pi(\x)\|\geq \delta}\exp(-\beta 
(F(\x)-M(\eta))d\x\to 0$ as $\beta\to +\infty$. The denominator 
\be 
\int_{\|\x-\pi(\x)\|\leq \eta}\exp\big(-\beta 
(F(\x)-M(\eta)\big)d\x>\int_{\|\x-\pi(\x)|\leq \eta}dx.
\ee 
It follows that  
\be 
P(\|X_\beta-\pi(X_\beta)\|\geq \delta)\leq 
\frac{\int_{\|\x-\pi(\x)\|\geq \delta}\exp\big(-\beta(F(\x)-M(\eta))\big)d\x}
{\int_{\|\x-\pi(\x)\|\leq \eta}dx}\to 0
\ee 
as $\beta\to +\infty$. \\

%\end{document} 
Now we are interested in the speed of convergence 
of $X_\beta-\pi(X_\beta)\to 0$ as $\beta\to +\infty$. 
The first step of this convergence is based on the following.  

{\bf Prop:}\label{Fsimplification} 
Let $\x(\y)$ be any Lasso estimator and $m=F(\x(\y))$ be the minimum of the objective function 
$F(\x)$ "equation~(\ref{F})". The function $F(\x)-m$ is equal to  
\ben\label{Fminusm} 
\sum_{i=1}^p|x_i|\big(1-sgn(x_i)\xi_i\big)+\frac{\|\A(\x-\x(\y))\|^2}{2}. 
\een 
And then 
\ben\label{partialI} 
\sum_{i=1}^p|x_i|\big(1-sgn(x_i)\xi_i\big)&=&
\sum_{i\in I} |x_i|\big(1-sgn(x_i)\xi_i\big)+\nonumber\\
&&2\sum_{i\in \partial I: sgn(x_i)\xi_i=-1} |x_i|. 
\een 
Here $\xi$ is defined by "equation~(\ref{xi})", $I$ and $\partial I$ are defined by "equation~(\ref{I})", and  "equation~(\ref{partialI})". 

{\bf Proof:} From the equality $\|\A\x-\y\|^2=\|\A(\x-\x(\y))\|^2+2\langle\A(\x-\x(\y)),\A\x(\y)-\y\rangle+
\|\A\x(\y,t)-\y\|^2$, we have  
\be
&&F(\x)=\\ 
&&\|\x\|_1+\frac{\|\A(\x-\x(\y))\|^2}{2}+\langle\A(\x-\x(\y)),\A\x(\y)-\y\rangle+\frac{\|\A\x(\y)-\y\|^2}{2}\\
&&=\|\x\|_1+\frac{\|\A(\x-\x(\y))\|^2}{2}+\langle\x-\x(\y),\A^*(\A\x(\y)-\y)\rangle+\frac{\|\A\x(\y)-\y\|^2}{2}.
\ee 
From the equality $\xi=\A^*(\y-\A\x(\y))$ , we have 
\ben\label{fc} 
\langle\x-\x(\y),\A^*(\A\x(\y)-\y)\rangle&=&-\langle\x-\x(\y),\xi\rangle\nonumber\\ 
&=&-\langle \x,\xi\rangle+\|\x(\y)\|_1.
\een 
Now formulas "equation~(\ref{Fminusm})" and "equation~(\ref{partialI})" are an easy consequence of the formula "equation~(\ref{fc})".

Now, we are interested in the asymptotic independence of 
the components $(X_\beta(i):\quad i\in I)$, 
$(X_\beta(i):\quad i\in \partial I)$ as $\beta\to +\infty$.  
We are going to solve this problem when 
$\A_{\partial I}^*\A_{\partial I}$ is invertible. 
In this case Lasso is a singleton $\{\x(\y)\}$. 

The support of $\x(\y)$ is $S=\{i:\quad x_i(\y)\neq 0\}$. 
The complementary of $S$ is $I_0=\{i: x_i(\y)=0\}$. The boundary of 
$\partial I_0=\{i: x_i(\y)=0,\,|\xi_i|=1\}$. The family
 $(S, I_0\setminus \partial I_0, \partial I_0)$ is a partition of $\{1, 2, \ldots, p\}$. In the sequel $\Rb^p$ is considered as the set of the sequences 
$(x_i: i\in (I_0\setminus \partial I_0)\cup \partial I_0\cup S)$   
indexed by 
$(I_0\setminus \partial I_0)\cup \partial I_0\cup S$. The notation 
$\Rb^J$ will denotes the set of the sequences $(x_j:\quad j\in J)$ 
with values in $\Rb$.   

Observe that 
$I=I_0\setminus\partial I_0$ "equation~(\ref{I})", and $\partial I=S\cup\partial I_0$ "equation~(\ref{partialI})". 
For $i\in S$ and for $x_i$ near $x_i(\y)$, we have $sgn(x_i)=\xi_i$. 
In this case the equality "equation~(\ref{partialI})" becomes 
\ben\label{l1lasso}
\sum_{i\in I_0\setminus \partial I_0} |x_i|\big(1-sgn(x_i)\xi_i\big)+
2\sum_{i\in\partial I_0: sgn(x_i)\xi_i=-1} |x_i|.
\een 

Now we decompose $X_\beta$ as following. Each partition $\partial I_0^-\cup\partial I_0^+$ of $\partial I_0$ defines two sets 
\be 
\Delta^-&:=&\Delta(\partial I_0^-)\\
&=&\big\{\x\in\Rb^p: \quad x_i\xi_i=-1, \forall\,i\in \partial I_0^-\big\},\\
\Delta^+&:=&\Delta(\partial I_0^+)\\
&=&\big\{\x\in\Rb^p: \quad x_i\xi_i=1, \forall\,i\in \partial I_0^+\big\}. 
\ee 
We have 
\be 
\Rb^p=\bigcup_{\partial I_0^-\cup\partial I_0^+=\partial I_0}\Delta^-\cap\Delta^+.
\ee  

It follows that for each suitable function $f$ 
\be 
\Eb[f(X_\beta)]=\sum_{\partial I_0^-\cup \partial I_0^+=\partial I_0}\Eb\big[f(X_\beta)\,|\,X_\beta\in \Delta^{-}\cap \Delta^+\big]
\Pb(X_\beta\in \Delta^{-}\cap \Delta^+).
\ee 
The main result of this section is the following. \\
{\bf prop:} \label{fondamentale} We have  
for each partition $K^-\cup K^+=\partial I_0$ with $K^-\neq \emptyset$ that 
\be 
\Pb(X_\beta\in \Delta(K^{-})\cap \Delta(K^+))\to 0\quad \mbox{as}\quad \beta\to +\infty.
\ee

{\bf proof:}  We suppose without loosing any generality for all $i\in\partial I_0$
that $\xi_i=1$. From "equation~(\ref{CVP})", 
we have for large $\beta$ that 
\be 
\Pb\big(X_\beta\in\Delta(K^{-})\cap \Delta(K^+) \big)\approx 
\frac{A(\beta,\delta,K^+,K^-)}{\sum_{\partial I_0^+\cup \partial I_0^-=\partial I_0} 
A_\beta(\delta,\partial I_0^+,\partial I_0^-)},
\ee 
where $\delta$ is small and 
\be 
A(\beta, \delta,\partial I_0^+,\partial I_0^-)&=&\int_{[\x\in 
\Delta^{-}\cap \Delta^+, \|\x-\x(\y)\|_\infty\leq\delta ]}
\exp\big(-\beta G(\x)\big)d\x,\\
G(\x)&=&\sum_{i\in I_0\setminus \partial I_0}|x_i|\big(1-\xi_i sgn(x_i)\big)+2\sum_{i\in\partial I_0^-}|x_i|+\\
&&\frac{\|\A(\x-\x(\y))\|^2}{2}.
\ee 
We recall that by hypothesis $K^-\neq \emptyset$, but in the denominator the sum 
$\sum_{\partial I_0^+\cup \partial I_0^-=\partial I_0}$ contains the case $\partial I_0^-= \emptyset$. 
 
We use the new variables 
\be 
u_i=\beta x_i,\quad i\in I_0\setminus \partial I_0^+,\\
v_i=\sqrt{\beta}(x_i-x_i(\y)),\quad i\in S\cup \partial I_0^+,
\ee 
and then we obtain 
\be 
A(\beta,\delta,\partial I_0^+,\partial I_0^-)=\beta^{-|I_0\setminus\partial I_0^+|-\frac{|S|+|\partial I_0^+|}{2}}
C(\beta,\delta, \partial I_0^+,\partial I_0^-),
\ee 
where 
\be 
C(\beta,\delta, \partial I_0^+,\partial I_0^-)=
\int_{[-\delta\beta, 0]^{\partial I_0^-}\times [-\delta\beta,\delta\beta]^{I_0\setminus \partial I_0}\times 
[-\delta\sqrt{\beta},\delta\sqrt{\beta}]^S\times 
[0,\delta\sqrt{\beta}]^{\partial I_0^+}}\\
\exp\left(-G(\u,\v,\beta,\partial I_0^+,\partial I_0^-)\right) d\u d\v, 
\ee 
with 
\be 
G(\u,\v,\beta,\partial I_0^+,\partial I_0^-)&&=
\sum_{i\in I_0\setminus \partial I_0}|u_i|\big(1-\xi_i sgn(u_i)\big)+2\sum_{i\in\partial I_0^-}|u_i|+\\
&&\frac{\|\A_{S\cup \partial I_0^+}\v_{S\cup \partial I_0^+}+\beta^{-1/2}
\A_{I_0\setminus \partial I_0^+}\u_{I_0\setminus\partial I_0^+}\|^2}{2}.
\ee 
Observe that $G(\u,\v,\beta,\partial I_0^+,\partial I_0^-)$ converges to 
\be 
G(\u,\v,\partial I_0^+,\partial I_0^-)=&&\sum_{i\in I_0\setminus \partial I_0}|u_i|\big(1-\xi_i sgn(u_i)\big)+\\
&&2\sum_{i\in\partial I_0^-}|u_i|+
\frac{\|\A_{S\cup \partial I_0^+}\v_{S\cup \partial I_0^+}\|^2}{2}
\ee 
as $\beta\to +\infty$, and then 
$C(\beta,\delta, \partial I_0^+,\partial I_0^-)$ converges to the following positive constant 
\be 
C(\partial I_0^+,\partial I_0^-):=
\int_{(-\infty, 0]^{\partial I_0^-}\times (-\infty,+\infty)^{(I_0\setminus \partial I_0)\cup S}\times(0,+\infty)^{\partial I_0^+}}\\
\exp\left(-G(\u,\v,\partial I_0^+,\partial I_0^-)\right) d\u d\v
\ee  
as $\beta\to +\infty$. By observing 
that $|\partial I_0^+|=|\partial I_0|$ is the minimizer of 
\be 
|\partial I_0^+|\to |I_0|-\frac{|\partial I_0^+|}{2}+\frac{|S|}{2}, 
\ee 
it follows that for $K^-\neq \emptyset$,
\be 
\frac{A(\beta,\delta,K^+,K^-)}{\sum_{\partial I_0^+\cup \partial I_0^-=\partial I_0} 
A(\beta,\delta,\partial I_0^+,\partial I_0^-)}
\ee 
converges to 0 as $\beta\to +\infty$. 

As a consequence we derive that as $\beta\to +\infty$, 
\be 
\Pb(X_\beta(i)\xi_i=1,\forall\,i\in \partial I_0)\to 1,
\ee  
and then we get the following. \\

{\bf Theo:}  \cite{DDN} 
If $\A_{\partial I}^*\A_{\partial I}$ is invertible, then the components  
\be 
\Big(\big(\beta X_\beta(i), i\in I_0\setminus\partial I_0\big), 
\big(\sqrt{\beta}(X_\beta(i)-x_i(\y)\big):\quad i\in S\cup \partial I_0\Big)
\ee 
are asymptotically independent as $\beta\to +\infty$. Their asymptotic PDF are 
proportional respectively to 
\be 
\prod_{i\in I_0\setminus\partial I_0} \exp\Big(-|x_i|(1-sgn(x_i)\xi_i)\Big),\\
\exp\Big(-\frac{\|\A_{S\cup \partial I_0}(\x-\x(\y))_{S\cup \partial I_0}\|^2}{2}\Big).
\ee

\section{Bayesian Lasso and multivalued diffusion}     
First we solve rigorously the following 
stochastic differential equation 
\ben 
d\x=-\big[\partial\|\x\|_1+\A^*(\A\x-\y)\big]dt+d\w, 
\label{mblasso}
\een 
where $\w$ is the standard Brownian motion. Second we show that 
the solution of "equation~(\ref{mblasso})" is ergodic with the stationary 
probability density  
"equation~(\ref{posterior})" with $\beta=1$.  

\subsection{Yosida approximation} 
Let $\varphi:\Rb^p\to (-\infty, +\infty]$ be a proper l.s.c. convex function, and  
$\mathcal{P}(\Rb^p)$ be the set of subsets of $\Rb^p$. 
The sub-differential $\partial\varphi$ is the map 
from $\Rb^p\to \mathcal{P}(\Rb^p)$ defined by 
\be 
\partial\varphi(\x)=\{\v\in\Rb^p:\quad \varphi(\x+\h)\geq \varphi(\x)+\langle h,\v\rangle,\,
\forall\,\h\in\Rb^p\}.
\ee  
The domain 
\be 
Dom(\partial\varphi)=\{\x:\quad\partial\varphi(\x)\neq\emptyset\}.
\ee 

A sequence of single valued approximations for the subdifferential 
$\partial\varphi(\x)$ is based on Yosida approximation. For 
each $\varepsilon >0$ and $\z\in\Rb^p$, the equation 
\be 
\x=\z+\varepsilon\partial\varphi(\z)
\ee 
has a unique solution denoted by  
\be 
\z&&=(\I+\varepsilon\partial\varphi)^{-1}(\x)\\
&&:=prox_{\varepsilon\varphi}(\x).
\ee 
The map $prox_{\varepsilon\varphi}: \Rb^p\to Dom(\partial\varphi)$ 
is called proximal function. 
The Yosida approximation of the sub-differential $\partial\varphi$ 
is the application 
\be 
\beta_{\varepsilon}(\x):=\frac{\x-prox_{\varepsilon\varphi}(\x)}{\varepsilon}. 
\ee 
The following are well known see e.g. \cite{Lepingle}.

{\bf prop:}  We have 
\begin{enumerate}
\item $prox_{\varepsilon\varphi}$ is a contraction from $\Rb^p$ to 
$Dom(\partial\varphi)$. 
\item $\beta_{\varepsilon}$ is monotone on the whole $\Rb^p$, i.e. 
\be 
\langle\beta_{\varepsilon}(\x_1)-\beta_{\varepsilon}(\x_2),\x_1-\x_2\rangle\geq 0,
\ee 
for all $\x_1, \x_2\in \Rb^p$, and is Lipschitz continuous with the constant $\frac{1}{\varepsilon}$.
\item For every $\x\in\Rb^p$, $\beta_{\varepsilon}(\x)\in\partial\varphi(
prox_{\varepsilon\varphi}(\x))$.  
 \end{enumerate}

{\bf prop:} For each $\varepsilon >0$, 
the map 
\be 
\x\in\Rb^p\to \varphi_{\varepsilon}(\x)=\min\{\varphi(\z)+\frac{\|\x-\z\|^2}{2\varepsilon}\}
\ee 
is called the Yosida approximation of the function $\varphi$. We have  
\begin{enumerate}
\item $\varphi_{\varepsilon}$ is convexe with the domain $\Rb^p$. 
\item $\varphi_{\varepsilon}$ is of class $C^1$ 
with $\nabla\varphi_{\varepsilon}=\beta_{\varepsilon}$. 
\item The infimum defining $\varphi_{\varepsilon}(\x)$ is attained at 
$prox_{\varepsilon\varphi}(\x)$, and 
\be 
\varphi_{\varepsilon}(\x)=\frac{\varepsilon}{2}\|\beta_\varepsilon(\x)\|^2+\varphi_{\varepsilon}(prox_{\varepsilon\varphi}(\x)).
\ee 
\item Letting $\varepsilon\downarrow 0$, we have $\varphi_{\varepsilon}\uparrow \varphi(\x)$ for all $\x\in\Rb^p$. 
\end{enumerate}

%\end{document} 
%We conclude this section by stating the implicit Euler scheme 
%see e.g. \cite{Lepingle}. 

%\begin{prop} Let $F$ be a convex and $C^1(\Rb^p,\Rb)$ function. The implicit algorithm scheme
%with the step size $\alpha >0$ 
%\be 
%\x(k+1)=\x(k)-\alpha\nabla\,F(\x(k+1))
%\ee 
%is defined by 
%\be 
%\x(k+1)=prox_{\alpha F}(\x(k)).
%\ee 
%Moreover, the sequence $(F(\x(k)):\,k=0, 1, \ldots)$
%is decreasing. 
%\end{prop} 
In the case $\varphi(\x)=\|\x\|_1$, 
we have 
\be 
prox_{\alpha\varphi}(\x)=(\x+\alpha){\bf 1}_{[\x\leq -\alpha]}+(\x-\alpha){\bf 1}_{[\x\geq\alpha]},
\ee 
and 
\be 
\varphi_{\varepsilon}(\x)&=&\min\Big\{\sum_{i=1}^p|z_i|+\frac{\|\z-\x\|^2}{2\varepsilon}\Big\}\\
&=&\sum_{i=1}^p\min\Big\{|z_i|+\frac{|z_i-x_i|^2}{2\varepsilon}\Big\}\\
&=&\sum_{i=1}^p \Big[(|x_i|-\frac{\varepsilon}{2}){\bf 1}_{[|x_i|\geq \varepsilon]}+
\frac{|x_i|^2}{2\varepsilon}{\bf 1}_{[|x_i|\leq \varepsilon]}\Big].  
\ee 
The gradient 
\be 
\nabla\varphi_{\varepsilon}(\x)&=&\beta_{\varepsilon}(\x)\\
&=&sgn(\x)
{\bf 1}_{[|\x|\geq \varepsilon]}+\frac{\x}{\varepsilon}{\bf 1}_{[|\x|\leq \varepsilon]}.
\ee 
Finally 
\be 
prox_{\alpha\varphi_{\varepsilon}}(\x)=(\x+\alpha){\bf 1}_{[\x\leq -\alpha-\varepsilon]}+\frac{\varepsilon \x}{\alpha+\varepsilon}{\bf 1}_{[|\x|\leq \alpha+\varepsilon]}+(\x-\alpha){\bf 1}_{[\x\geq \alpha+\varepsilon]}.
\ee  
%\end{document}
\subsection{Multivalued stochastic differential equation} 
Now, we come back to Multivalued stochastic differential equation. 
Let $\w$ be the standard Brownian motion on $\Rb^p$ and 
$\b:\Rb^p\to \Rb^p$ be a smooth map. 
A solution of the $\Rb^p$-multivalued stochastic differential equation (abbreviated MSDE) 
\ben 
d\x_t=-\partial\varphi(\x_t)dt-\b(\x_t)dt+d\w_t
\label{MSDE} 
\een 
is a couple of continuous adapted stochastic processes  $t\in [0,+\infty)\to (\x(t),\l(t))$ 
with values in $\Rb^p\times\Rb^p$, and such that $\l(0)=0$, 
$t\to \l(t)$ has bounded variation on each compact interval and 
\be 
d\x_t&=&-d\l_t-\b(\x_t)dt+d\w_t,
\ee 
and "$\frac{d\l(t)}{dt}\in\partial\varphi(\x(t))$", i.e. 
the measure $\langle\x_t-\alpha_t,d\l_t-\beta_tdt\rangle$ is non-negative  
for all continuous trajectory $t\to (\alpha_t,\beta_t)$ 
such that $\beta_t\in\partial\varphi(\alpha_t)$. Observe that if $d\l_t=\l_t'dt$, 
then $\l_t'\in\partial\varphi(\x_t)$.  

It's known that if 
\be 
\|\b(\x_1)-\b(\x_2)\|&&\leq C\|\x_1-\x_2\|,\quad \forall\,\x_1, \x_2,\\
\|\b(\x)\|&&\leq C(1+ \|\x\|),\quad \forall\,\x, 
\ee 
then there exits a unique solution $(\x,\l)$. See e.g. \cite{CEPAthese},\cite{CEPAAnal}, 
\cite{Bernardin}, \cite{Kree}, \cite{Bensoussan}, 
\cite{Storm}. It follows that "equation~(\ref{mblasso})" has a unique solution 
$(\x,\l)$. In general the measure $d\l_t$ is not absolutely continuous with respect to the Lebesgue measure $dt$. 
However we are going to show that $d\l_t$ is absolutely continuous in the case 
"equation~(\ref{mblasso})". We recall two methods for constructing the solution $\x$ of "equation~(\ref{mblasso})". 

1) By choosing $\varphi(\x)=\|\x\|_1$, $\b(\x)=\A^*(\A\x-\y)$, then 
the solution of "equation~(\ref{mblasso})" is the unique couple $(\x,\l)$ of continuous maps  such that $\l(0)=0$, $t\to \l(t)$ has bounded variation on each compact interval and 
\ben 
d\x(t)=-[d\l_t+\A^*(\A\x(t)-\y)dt]+d\w_t ,\quad \frac{d\l(t)}{dt}\in\partial\|\x(t)\|_1. 
\label{lskorokhod} 
\een 

2) By choosing $\varphi(\x)=\|\x\|_1+\frac{\|\A\x-\y\|^2}{2}$, 
$\b(\x)=0$, then the solution of "equation~(\ref{mblasso})"
is given by the couple $(\x(t),\k(t))$ such that 
\be 
d\x(t)=-d\k(t)+d\w(t),\quad \k(t)\in\partial\varphi(\x(t)). 
\ee 

The uniqueness of the solution of "equation~(\ref{mblasso})" 
implies that $d\k(t)=d\l_t+\A^*(\A\x(t)-\y)dt$. 
Now, we are going to show that 
$\l$ is absolutely continuous. For this aim we recall Skorokhod problem 
\cite{CEPAAnal}. Let $\f$ be any continuous function from $[0,T]\to \Rb^d$, 
and $\psi:\Rb^p\to \Rb$ be any convex function.  
Then there exists 
a unique couple $(\x,\k)$ of continuous maps  such that $\k(0)=0$, 
$t\to \k(t)$ has bounded variation on each compact interval, 
\ben 
\x(t)=\f(t)-\k(t),\quad \forall\,t\geq 0,
\label{skorokhod} 
\een 
and the measure $\langle\x(t)-\alpha(t),d\k(t)-\beta(t)dt$ is nonnegative  
for all continuous trajectory $t\to (\alpha(t),\beta(t))$ 
such that $\beta(t)\in\partial\psi(\alpha(t))$. Now we are ready to announce our result. 
{\bf prop:} Suppose that 
\ben 
m=\sup\{\|\v\|:\quad \v\in\bigcup_{\x\in\Rb^p}\partial\psi(\x)\}
\label{m} 
\een 
is finite. Then the function $\l$ solution of Skorokhod problem "equation~(\ref{skorokhod})" is 
absolutely continuous. \\

{\bf proof:} Let $\e\in\Rb^p$ such that $\|\e\|=1$, $\gamma >0$ and 
$v\in\partial\psi(\gamma\e)$ having the smallest Euclidean norm. 
As $(\x,\k)$ is the solution of Skorokhod problem, then 
\be 
\langle\x(t)-\gamma\e,d\l(t)\rangle&\geq& \langle\x(t)-\gamma\e,\v dt\rangle\\
&\geq & -m\left(\|\x(t)\|+\gamma\right)dt.
\ee 
For each $0\leq s < t$, we have 
\be 
\langle \l(t)-\l(s), \e\rangle&=&\int_{s}^{t}\langle\e,d\l(u)\rangle\\
&=&\gamma^{-1}\int_{s}^{t}\langle\x(u),d\l(u)\rangle-\gamma^{-1}
\int_{s}^{t}\langle\x(u)-\gamma\e,d\l(u)\rangle\\
&\leq & \gamma^{-1}\int_{s}^{t}\langle\x(u),d\l(u)\rangle+m\gamma^{-1}\int_{t_i}^{t_{i+1}}\|\x(u)\|du+m (t-s).
\ee 
From the latter inequality and 
\be 
\|\l(t)-\l(s)\|=\sup\{\langle \l(t)-\l(s), \e\rangle:\quad \e\in\Rb^p,\quad \|\e\|=1\},
\ee 
and by tending $\gamma\to +\infty$, we get  
\be 
\|\l(t)-\l(s)\|\leq m (t-s). 
\ee 
Which achieves the proof. \\

By choosing $\f(t)=\x_0-\int_0^t \A^*(\A\x(s)-\y)ds]+\w_t$, 
we derive that $(\x,\l)$ "equation~(\ref{lskorokhod})" is the solution of Skorokhod problem. As the hypothesis "equation~(\ref{m})" is satisfied for $\psi(\x)=\|\x\|_1$,
with $m=1$,  
then $\l$ is absolutely continuous. Finally the solution 
of "equation~(\ref{mblasso})" satisfies 
\ben 
\x(t)=\x(0)-\int_0^t\big[\v(s)+\A^*(\A\x(s)-\y)\big]ds+\w(t), 
\label{xv} 
\een 
and $\v(t)\in\partial\|\x(t)\|_1$, $\|\v(t)\|\leq 1$, $dt$ a.e.
Moreover we can show that a.s. for $i=1, \ldots, p$ that  
$x_i(t)\neq 0$ and $v_i(t)=sgn(x_i(t))$, $dt$ a.e. 
The "equation~(\ref{xv})" becomes 
\ben 
d\x(t)=\frac{1}{2}\nabla\ln\big(\rho(x(t))\big)dt+d\w(t),  
\label{sgnx} 
\een 
where 
\ben 
\rho(\x):=\frac{1}{Z}\exp\big(-2\|\x\|_1-\|\A\x-\y\|^2\big).
\label{invariantdensity} 
\een 
The equation "equation~(\ref{sgnx})" is known as distorted Brownian motion \cite{Fukushima} with the 
generalized Schr\"odinger operator 
\be 
H=-\frac{1}{2}\Delta-\frac{1}{2}
\Big(\sum_{i=1}^p\delta(x_i)+Trace(\A^*\A)\Big)+\frac{1}{2}\|sgn(\x)+\A^*(\A\x-\y)\|^2.
\ee 
Here $\Delta$ is Laplacian operator and $\delta$ denotes the Dirac measure at 0. 
\subsection{Transition probabilities in the one dimensional case} 
In the one dimensional case 
\be 
dx(t)=-\lambda sgn(x(t))dt+dw(t), \quad x(0)=x_0,\quad \lambda >0 
\ee 
is known as bang-bang Brownian motion \cite{Shreve}, 
or the diffusion with V potential \cite{Risken}. In this case 
Schr\"odinger operator has the form 
\be 
H=-\frac{1}{2}\frac{d^2}{d^2x}+\frac{1}{2}(1-\delta).
\ee 
The transition probabilities $p^\lambda(x,t\,|\,x_0,0)$ of the bang-bang Brownian motion 
is known \cite{Borodin}. We can calculate it using Girsanov Formula,  
and the trivariate density of Brownian motion, its local time and occupation times
(\cite{Karatzas}). We obtain  
\be 
p^\lambda(x,t\,|\,x_0)=q^{\lambda}(x,t\,|\,x_0)\lambda\exp(-2\lambda |x|)
\ee 
where 
\be 
q^{\lambda}(x,t\,|\,x_0)&=&\exp\Big(\lambda(|x_0|+|x|)-\frac{t\lambda^2}{2}\Big)\gamma_t(x-x_0)+F\Big(\frac{\lambda t-(|x|+|x_0|)}{\sqrt{t}}\Big),\\ 
F(x)&=&\int_{-\infty}^x\frac{exp(-\frac{u^2}{2})}{\sqrt{2\pi}}du,\\
\gamma_t(u)&=&\frac{exp(-\frac{u^2}{2t})}{\sqrt{2t\pi}}. 
\ee  
Observe that $p^\lambda(x,t\,|\,x_0,0)\to \lambda\exp(-2\lambda |x|)$ 
as $t\to +\infty$ for all $x_0$. Hence, the MSDE 
\be 
dx(t)=-\lambda sgn(x(t))dt+dw(t)
\ee 
is ergodic with the invariant density $\lambda\exp(-2\lambda|x|)$. 

\section{Sampling using multivalued SDE} 
As we said before, the solution $(\x(t))$
of "equation~(\ref{mblasso})" is ergodic. 
It follows that $\lim_{T\to +\infty}\x(T)$
has the probability distribution $\rho$ "equation~(\ref{invariantdensity})".  
If we dispose of a trajectory $t\in [0,T]\to\x_t$ for large $T$, 
then for any 
$\rho$-integrable function $h$,  
\be 
\frac{1}{T}\int_{0}^Th(\x_t)dt\approx \Eb[h(\x(T))]\approx \int_{\R^p} h(\x)\rho(\x)d\x.
\ee 
Hence for large $T$ the expectation $\Eb[\x(T)]$ of 
the solution "equation~(\ref{mblasso})" is close to 
Bayes estimator "equation~(\ref{bayes})". 
 We will approximate $\Eb[\x(T)]$ using numerical schemes of 
"equation~(\ref{mblasso})" and the timestep 
\ben 
\Delta t_l=2^{-l}T, 
\label{timestep} 
\een 
with the level  $l=l_s, l_s+1, \ldots$. In all the sequel 
the small level $l_s:=\frac{\ln(T)}{\ln(2)}+1$.  

Having a numerical scheme $(\x_L(sc,k):\quad k=1, \ldots, 2^L)$
such that $\Eb[\x_L(sc,2^L)]\to \Eb[\x(T)]$ as $L\to +\infty$,  
we need to calculate $\Eb[\x_L(sc,2^L)]$ for large $L$. 
To achieve this goal we use 
Monte Carlo (MC) and multilevel Monte Carlo (MLMC) algorithms.
We will discuss the efficiency of MC and MLMC estimates.  
We will mimic the results obtained  
in \cite{Rosin} for Coulomb collisions, and propose 
a method for calculating the cost.    
\section{MC Efficiency and computational cost} 
Given a sample $(\x_l^k(sc,2^l):\,k=1, \ldots, N_l)$
of $\x_l(sc,2^l)$ having the size $N_l$, we define
\ben 
&&\hat{\x}_l^{N_l}(sc,2^l)=\frac{1}{N_l}\sum_{k=1}^{N_l}\x_l^k(sc,2^l),\quad l,\,\mbox{and}\, N_l\,\,
\mbox{are fixed},\label{moyenneempiriquedexl}\\
&&\hat{\x}_l(sc,2^l):=\Eb[\x_{l}(sc,2^l)]=\lim_{N_l\to +\infty}\hat{\x}_l^{N_l}(sc,2^l),
\label{moyennedexl}\\ 
&&\hat{\x}(T):=\Eb[\x(T)]=\lim_{l\to +\infty}
\hat{\x}_l(sc,2^l). 
\label{moyennedexT}
\een 
We recall that MC proposes to estimate  $\hat{\x}_l(sc,2^l)$ 
"equation~(\ref{moyennedexl})"
by $\hat{\x}_l^{N_l}(sc,2^l)$  "equation~(\ref{moyenneempiriquedexl})". 

If we estimate $\hat{\x}(T)$ "equation~(\ref{moyennedexT})" by $\hat{\x}_l^{N_l}(sc,2^l)$, then the error has two sources.  
The approximation of $\Eb[\x(T)]$ by $\Eb[\x_l(sc,2^l)]$, 
and a finite sampling error that depends on the number of samples $N_l$. 

An accurate estimate $\hat{\x}_l^{N_l}(sc,2^l)$ of 
$\hat{\x}(T)$ is one for which the mean square error 
\be 
&&MSE=\Eb\Big[\|\hat{\x}(T)-\hat{\x}_l^{N_l}(sc,2^l)\|^2\Big]\\
&&=\|\hat{\x}(T)-\hat{\x}_l(sc,2^l)\|^2+
\Eb\Big[\|\hat{\x}_l(sc,2^l)-\hat{\x}_l^{N_l}(sc,2^l)\|^2\Big]
\ee 
is small. We have 
\be 
\Eb\Big[\|\hat{\x}_l(sc,2^l)-\hat{\x}_l^{N_l}(sc,2^l)\|^2\Big]:=\frac{Var_l(sc)}{N_l},
\ee
where 
\be 
Var_l(sc)=\sum_{i=1}^pVar(x_{l,i}(sc,2^l)).
\ee 
Here $x_{l,i}(sc,2^l)$ is the $i$-th component of $\x_l(sc,2^l)$ and 
$Var(x_{l,i}(sc,2^l))$ its variance. 

The quantity 
\ben 
\|\hat{\x}(T)-\hat{\x}_l(sc,2^l)\|^2=e(sc,\Delta t_l)
\label{el} 
\een 
is a function of the timestep $\Delta t_l$.  
It is central in the computational cost and we suppose that 
is known.    

%The quantity $V_l(s)$ is estimated by 
%\be 
%\frac{1}{N}\sum_{k=1}^{N}\|\hat{\x}^{N}(2^l,s)-\x^k(2^l,s)\|^2
%\ee 
%with a large $N$.  
The estimate $\hat{\x}_l^{N_l}(sc,2^l)$ is accurate 
to within $\eta^2$ of $\hat{\x}(T)$ if 
\ben 
MSE&=&\Eb\Big[\|\hat{\x}(T)-\hat{\x}_l^{N_l}(sc,2^l)\|^2\Big]=\eta^2\nonumber\\
&=&e(sc,\Delta t_l)+\frac{Var_l(sc)}{N_l}. 
\label{constrainteta}
\een
The computational cost $K$ of obtaining $(\x_l^k(sc,2^l):\quad k=1, \ldots, N_l)$ is the product of the number of timestep 
$\frac{T}{\Delta t_l}=2^l$ and the number of samples $N_l$. Namely,
\be 
K(N_l,\Delta t_l)=N_l\frac{T}{\Delta t_l}=N_l2^l. 
\ee 
To make the scheme as efficient as possible, $K$ must be minimal subject
to the constraint "equation~(\ref{constrainteta})".  
Applying the method of Lagrange multipliers 
\be 
L(N_l,\Delta t_l,\lambda)=N_l\frac{T}{\Delta t_l}+\lambda\Big(e(sc,\Delta t_l)+\frac{Var_l(sc)}{N_l}-\eta^2\Big),
\ee 
we get the optimal choice 
\ben 
\frac{T}{\Delta t_l}-\lambda\frac{Var_l(sc)}{N_l^2}=0,\nonumber\\
-N_l\frac{T}{(\Delta t_l)^2}+\lambda\frac{\partial e}{\Delta t_l}
(sc,\Delta t_l)=0,\nonumber\\
e(sc,\Delta t_l)+\frac{Var_l(sc)}{N_l}=\eta^2.\label{Nlversiongeneral}
\een
It follows that 
\ben 
\frac{\partial\,e(sc,\Delta t_l)}{\Delta t_l}&=&\frac{\eta^2-e(sc,\Delta t_l)}{\Delta t_l},\label{deltae}\\
e(sc,\Delta t_l)&<& \eta^2. 
\een 
We propose to solve the latter system numerically as follows. 
In all the sequel we estimate $\hat{\x}(T)$ by $\hat{\x}_L(sc,2^L)$ 
with $L=16$. 
Hence we obtain the following approximation:   
\ben 
\|\hat{\x}_{L}(sc,2^L)-\hat{\x}_{l}(sc,2^l)\|^2\approx e(sc,\Delta t_l).\label{eLl} 
\een 
Second  
\be 
\frac{\partial\,e(sc,\Delta t_l)}{\Delta t_l}\approx \frac{e(sc,\Delta t_{l})-e(sc,\Delta t_{l+1})}{T2^{-l-1}}.
\ee 
The "equation~(\ref{deltae})" becomes 
\be 
3e(sc,\Delta t_{l})-2e(sc,\Delta t_{l+1})\approx \eta^2. 
\ee
Now we calculate for $l\geq l_s$ the quantity 
\ben 
3e(sc,\Delta t_l)-2e(sc,\Delta t_{l+1})
\label{3e-2e} 
\een 
until it becomes close to $\eta^2$ and 
\ben 
e(sc,\Delta t_l)< \eta^2.
\label{eleta}
\een 
Having $l$, we calculate 
$Var_l(sc)$ by 
\ben 
\sum_{i=1}^p\frac{1}{N}\sum_{k=1}^N\Big|x_{l,i}^k(sc,2^l)-\frac{1}{N}
\sum_{k=1}^Nx_{l,i}^k(sc,2^l)\Big|^2. 
\label{varlsc} 
\een 
Having $l$ and $Var_l(sc)$ we calculate the optimal 
sample size $N_l$ using the "equation~(\ref{Nlversiongeneral})" and then 
we derive the optimal cost 
$K_l$.   

\section{MLMC Efficiency and computational cost} 

Multilevel Monte Carlo (MLMC) was initially developed for financial mathematics \cite{Giles1}, \cite{Giles2} and now used in a disparate areas.  

Multilevel Monte Carlo considers multilevels. 
In our study we consider the levels $l=l_s, l_s+1, \ldots, l_m< L=16$. The smallest level $l_s$ is  choosen 
such that $\Delta t_{l_s}=\frac{1}{2}$.  
We generate a sample $(\x_{l_s}^k(sc,2^{l_s}): k=1, \ldots, N_{l_s})$ of size 
$N_{l_s}$ of $\x_{l_s}(sc,2^{l_s})$, and for each $l=l_s+1, \ldots, l_m$, 
we generate from the same underlying stochastic path and initial conditions
the samples $(\x_l^k(sc,2^l): k=1, \ldots, N_l)$ 
and $(\x_{l-1}^k(sc,2^{l-1}): k=1, \ldots, N_{l})$ 
respectively of $\x_l(sc,2^l)$ and  $\x_{l-1}(sc,2^{l-1})$. Moreover, the samples 
$(\x_{l_s}^k(sc,2^{l_s}): k=1, \ldots, N_{l_s})$, $(\x_l^k(sc,2^l), \x_{l-1}^k(sc,2^{l-1}): k=1, \ldots, N_{l})$ 
for $l=l_s+1$, \ldots, $l_m$ have to be independent.   
Using the telescoping sum 
\be 
\hat{\x}_{l_m}(sc,2^{l_m})=\hat{\x}_{l_s}(sc,2^{l_s})+\sum_{l=l_s+1}^{l_m}(\hat{\x}_l(sc,2^l)-
\hat{\x}_{l-1}(sc,2^{l-1})), 
\ee 
MLMC proposes the estimate    
\be 
\hat{\x}_{l_m}^{N_{l_m}}(2^{l_m})=\hat{\x}_{l_s}^{N_{l_s}}(sc,2^{l_s})+\sum_{l=l_s+1}^{l_m}\left(\hat{\x}_{l}^{N_l}(sc,2^l)-\hat{\x}_{l-1}^{N_{l}}(sc,2^{l-1})\right) 
\ee
of $\hat{\x}_{l_m}(sc,2^{l_m}):=\Eb[\x_{l_m}(sc,2^{l_m})]$.   

We introduce for each level $l$ and sample size $N_l$ the following notations: 
\be 
\hat{\x}_{l}^{N_l}(sc,2^l)-\hat{\x}_{l-1}^{N_{l}}(sc,2^{l-1}):=\delta\hat{\x}_{l}^{N_l}(sc,2^l). 
\ee 
It follows that
\ben 
\hat{\x}_{l_m}^{N_{l_m}}(sc,2^{l_m})&=&\hat{\x}_{l_s}^{N_{l_s}}(sc,2^{l_s})+\sum_{l=l_s+1}^{l_m}\left(\hat{\x}_{l}^{N_l}(sc,2^l)-\hat{\x}_{l-1}^{N_{l}}(sc,2^{l-1})\right)\nonumber\\
&:=&\hat{\x}_{l_s}^{N_{l_s}}(sc,2^{l_s})+\sum_{l=l_s+1}^{l_m}
\delta\hat{\x}_{l}^{N_l}(sc,2^l),
\label{MLMC} 
\een 
where $\delta\hat{\x}_{l}^{N_l}(sc,2^l):=\hat{\x}_{l}^{N_l}(sc,2^l)-\hat{\x}_{l-1}^{N_{l}}(sc,2^{l-1})$. 

An accurate estimate $\hat{\x}_{l_m}^{N_{l_m}}(sc,2^{l_m})$ of 
$\hat{\x}(T)$ is one for which the mean square error 
\be 
MSE:=\Eb\Big[\|\hat{\x}(T)-\hat{\x}_{l_m}^{N_{l_m}}(sc,2^{l_m})\|^2\Big]=
\|\hat{\x}(T)-\hat{\x}_{l_m}(sc,2^{l_m})\|^2+Var(\hat{\x}_{l_m}^{N_{l_m}}(sc,2^{l_m}))
\ee 
is small. 

If we set $V_{l_s}=Var(\x_{l_s}(sc,2^{l_s}))$, and  
for $l=l_s+1, \ldots, l_m$, 
\ben 
V_l=Var(\delta\x_{l}(sc,2^l)), 
\label{vl}
\een 
then 
\be 
Var(\hat{\x}_{l_m}^{N_{l_m}}(sc,2^{l_m}))=\sum_{l=l_s}^{l_m}\frac{V_l}{N_l},
 \ee 
and 
\be 
MSE&=&\|\hat{\x}(T)-\hat{\x}_{l_m}(sc,2^{l_m})\|^2+
\sum_{l=l_s}^{l_m}\frac{V_l}{N_l}. 
\ee 
If  
\ben 
\|\hat{\x}(T)-\hat{\x}_{l_m}(sc,2^{l_m})\|^2:=e(sc,\Delta t_{l_m})+
\sum_{l=l_s}^{l_m}\frac{V_l}{N_l}=\eta^2,\label{MLMCconstraint} 
\een 
then efficiency of MLMC is equivalent to minimize 
\ben 
K=\sum_{l=l_s}^{l_m}K_l=\sum_{l=l_s}^{l_m}N_l\frac{T}{\Delta t_l},
\label{KMLMC} 
\een 
under the constraint "equation~(\ref{MLMCconstraint})".

We estimate for $l_s\leq l< L=16$, $\|\hat{\x}(T)-\hat{\x}_{l}(sc,2^{l})\|^2$ by
$\|\hat{\x}_L(sc,2^L)-\hat{\x}_{l}(sc,2^{l})\|^2$ 
and then we are interested in the set $l(\eta)$ 
of levels $l$ such that   
\be 
e(sc,\Delta t_l)=\|\hat{\x}_L(sc,2^L)-\hat{\x}_{l}(sc,2^{l})\|^2 \approx \frac{\eta^2}{2}.   
\ee 
For each $lopt\in l(\eta)$, the "equation~(\ref{MLMCconstraint})" 
becomes 
\ben 
\sum_{l=l_s}^{lopt}\frac{V_l}{N_l}=\eta^2-e(sc,\Delta t_{lopt}). \label{Nopt}
\een 
Having the optimal $lopt$, the minimization of $K$ "equation~(\ref{KMLMC})" under the constraint (\ref{Nopt}) 
is solved by Lagrange multiplier 
\be 
\partial_{N_l}(K+\lambda(\frac{V_{l_s}}{N_{l_s}}+\sum_{l=l_s+1}^{lopt}\frac{V_l}{N_l}-
(\eta^2-e(sc,\Delta t_{lopt})))=0,\quad l=l_s, \ldots, lopt. 
\ee 
Hence 
\be 
2^l&=&\lambda \frac{V_l}{N_l^2},\quad l=l_s, \ldots, lopt,\\
\sum_{l=l_s+1}^{lopt}\frac{V_l}{N_l}&=&\eta^2-e(sc,\Delta t_{lopt}).
\ee 
It follows  for $l=l_s, \ldots, lopt$, that $N_l=\sqrt{\lambda V_l2^{-l}}$, 
and then  
\be 
\sum_{l=l_s}^{lopt}\frac{\sqrt{V_l2^l}}{\sqrt{\lambda}}=\eta^2-e(sc,\Delta t_{lopt}).
\ee 
Having $lopt$, we estimate $V_{l_s}$, and $(V_l: l=l_s+1, \ldots, lopt)$ by 
\ben 
\hat{V}_{l_s}=\sum_{i=1}^p
\frac{1}{N}\sum_{k=1}^N|x_{l_s,i}^k-\hat{x}_{l_s,i}^k|^2,\\ 
\hat{V}_l:=\sum_{i=1}^p
\frac{1}{N}\sum_{k=1}^N|\delta x_{l,i}^k-\delta\hat{x}_{l,i}|^2,\quad 
l=l_s+1, \ldots, lopt. 
\label{VL} 
\een 
Hence for $l=l_s, \ldots, lopt$ 
\ben 
N_l=\frac{1}{\eta^2-e(sc,\Delta t_{lopt})}\sqrt{V_l2^{-l}}\sum_{k=l_s}^{l}\sqrt{V_k2^k}.
\label{optimalNlL} 
\een 

Now, we are going to present our schemes.  
 
\section{Semi-implicit Euler schemes}
Numerical approximation has been tackled in \cite{Bernardin}, \cite{Lepingle}, \cite{Asiminoaei}, \cite{Pettersson}, see 
also \cite{Talay1, Talay2}. 
Semi-implicit Euler scheme (SIES) of "equation~(\ref{MSDE})"  
is given by  
\be 
\x_{l}(k+1)-\x_{l}(k)=-\nabla\varphi(\x_l(k+1))\Delta t_l-b(\x_{l}(k))\Delta t_l+\sqrt{\Delta t_l}\n(k+1), 
\ee
where $(\n(k+1):\quad k=0, 1, \ldots, 2^l-1)$ is a sequence of i.i.d. standard Gaussian vectors. 
Known $\x_{l}(k)$ and 
$\n(k+1)$, we have   
\ben 
\x_{l}(k+1)=prox_{\Delta t_l\varphi}\Big(\x_l(k)-
b(\x_{l}(k))\Delta t_l+\sqrt{\Delta t_l}\n(k+1)\Big). 
\label{SIES}
\een 

The weak and the strong convergence propertie of the scheme "equation~(\ref{SIES})" 
to the solution "equation~(\ref{MSDE})" 
are defined respectively in terms of 
\ben
e_w(\Delta t_l)=\|\Eb[\x(T)-\x_{l}(2^l)]\|,\label{sies2weakconvergence}\\
e_s(\Delta t_l)=\Eb\Big[\|\x(T)-\x_{l}(2^l)\|^2\Big]^\frac{1}{2}.
\label{sies2strongconvergence} 
\een 
From (\cite{Bernardin}) the strong error "equation~(\ref{sies2strongconvergence})" is estimated by 
\ben 
O((\Delta t_l\ln(\frac{1}{\Delta t_l})^\frac{1}{4}). 
\label{Bernardinestimate} 
\een 
By setting $\varphi(\x)=\|\x\|_1$, $\b(\x)=\A^*(\A\x-\y)$, 
the scheme "equation~(\ref{SIES})" is known as STMALA algorithm (\cite{Fort}). 
\section{Explicit Euler scheme}
\subsection{Algorithm EES1} 
By setting $\varphi(\x)=\|\x\|_1+\frac{\|\A\x-\y\|^2}{2}$, 
and $\b=0$, EES1 of (\ref{mblasso}) is given by  
\be 
\x_l(k+1)=prox_{\Delta t_l\varphi}(\x_l(k))+\sqrt{\Delta t_l}\n_{k+1}.
\ee  
The proximal $prox_{\Delta t_l\varphi}(\x^{(k)})$ is not computable, but 
for large $l$, we have 
\be 
prox_{\Delta t_l\varphi}(\x)\approx prox_{\Delta t_l\|\cdot\|_1}\Big(\x+\A^*(\y-\A\x)\Delta t_l\Big). 
\ee 
Finally we get
\ben 
\x_l(k+1)=prox_{\Delta t_l\|\cdot\|_1}\Big(\x_l(k)+\A^*\big(\y-\A\x_l(k)\big)\Delta t_l\Big)+\sqrt{\Delta t_l}\n_{k+1}, 
\label{PULA} 
\een 
known as PULA algorithm \cite{Pereyra}. 
\subsection{Algorithm EES2} 
By setting $\varphi(\x)=\|\x\|_1$, 
and $\b(\x)=\A^*(\A\x-\y)$, 
we obtain our new scheme 
\ben 
\x_l(k+1)=prox_{\Delta t_l\|\cdot\|_1}\big(\x_l(k)\big)+
\A^*\big(\y-\A\x_l(k)\big)\Delta t_l+\sqrt{\Delta t_l}\n_{k+1}.
\label{EE2} 
\een  

\section{Numerical implementation}

As an illustration we consider the case $p=10$, $n=7$ and the entries of the matrix $\A$ are independent Bernoulli random variables with values $\pm\frac{1}{\sqrt{n}}$, and $\w  \sim \mathcal{N}(0,\frac{1}{2}\I_n)$. We simulate the vector $\x(true)$ from 
the PDF $\exp(-2\|\x\|_1)$. We get the data $\y:=\A\x(true)+\w$ from  
a realization of $\A$ and $\w$. The time horizon $T=10$ the maximal level $L=16$ and the smallest level  
$l_s=5$.  
\subsection{Graphics of Trajectories of each scheme} 
For each scheme $sc$ and for each level $l=l_s, l_s+1, l_s+2$, 
we  plot the trajectories $k\in [0, 2^l]\to \x_l(sc,k)$.
For the largest level $L=16$ we plot only the first component. 
   
 %%%%%%%%%%%%%%%%%%%%%%%%%%%%%%%%%%%%%%%%%%%%%%%%%%%%%%%%%%%%%
\begin{figure}[!ht]
  \centering

     \includegraphics[width=16cm]{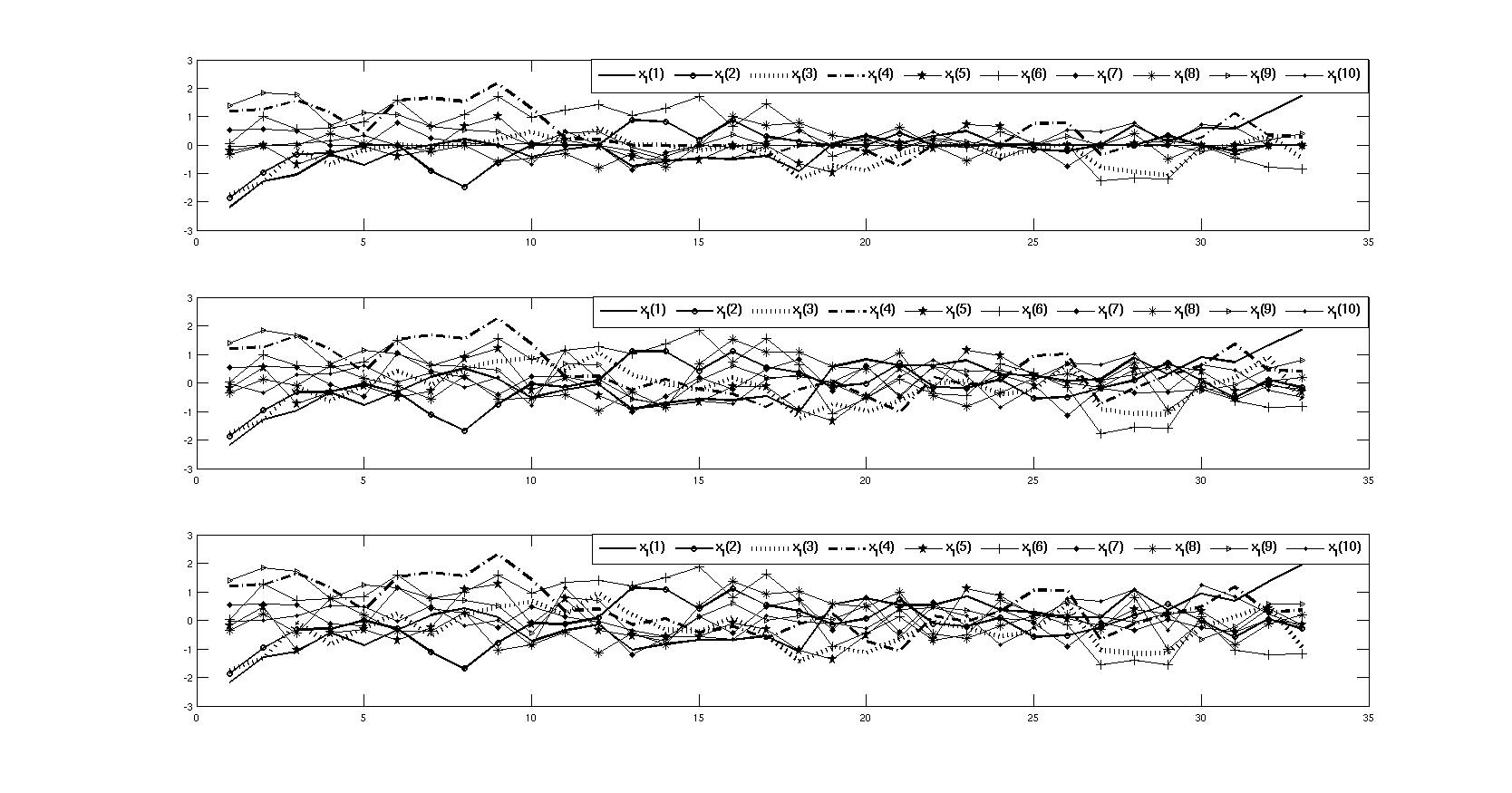}
 \caption{The chains of SIES, EES1 and EES2  for $l=l_s$.}

\end{figure}
\begin{figure}[!ht]
  \centering

     \includegraphics[width=16cm]{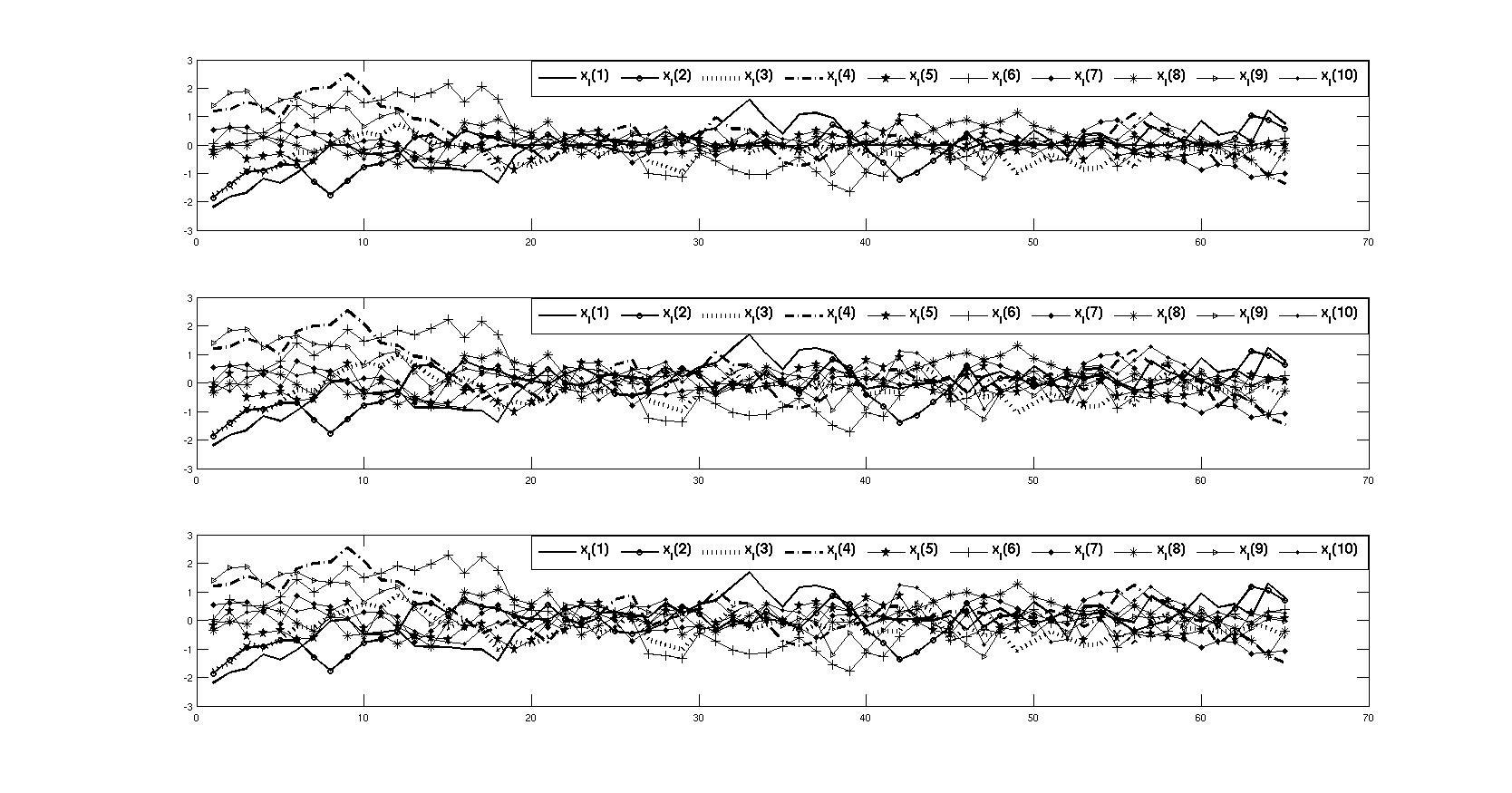}
 \caption{The chains of SIES, EES1 and EES2  for $l=l_s+1$.}

\end{figure}
\begin{figure}[!ht]
  \centering

     \includegraphics[width=13.7cm]{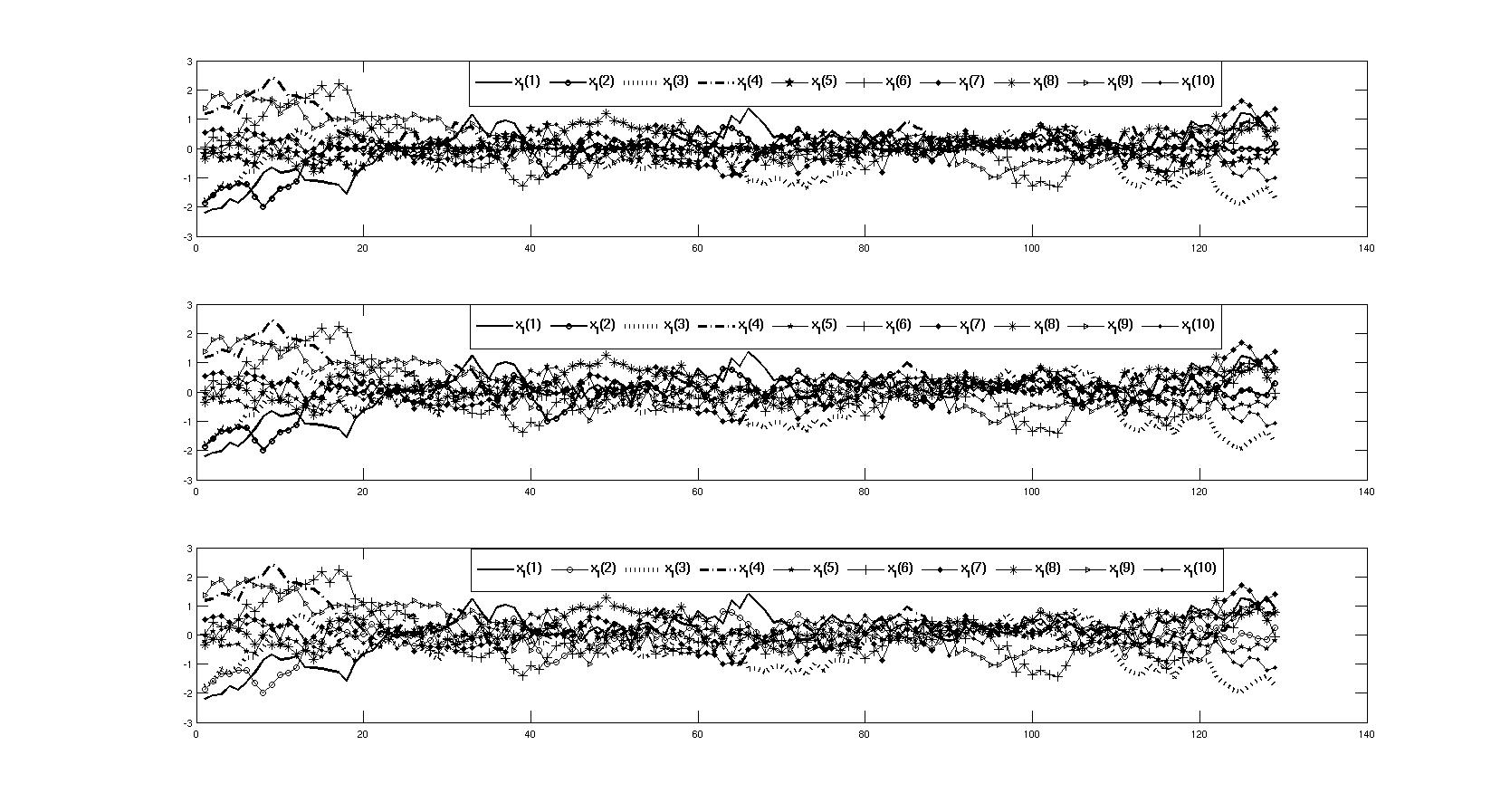}
 \caption{The chains of SIES, EES1 and EES2  for $l=l_s+2$.}
\end{figure}

\begin{figure}[!ht]
  \centering

     \includegraphics[width=13.7cm]{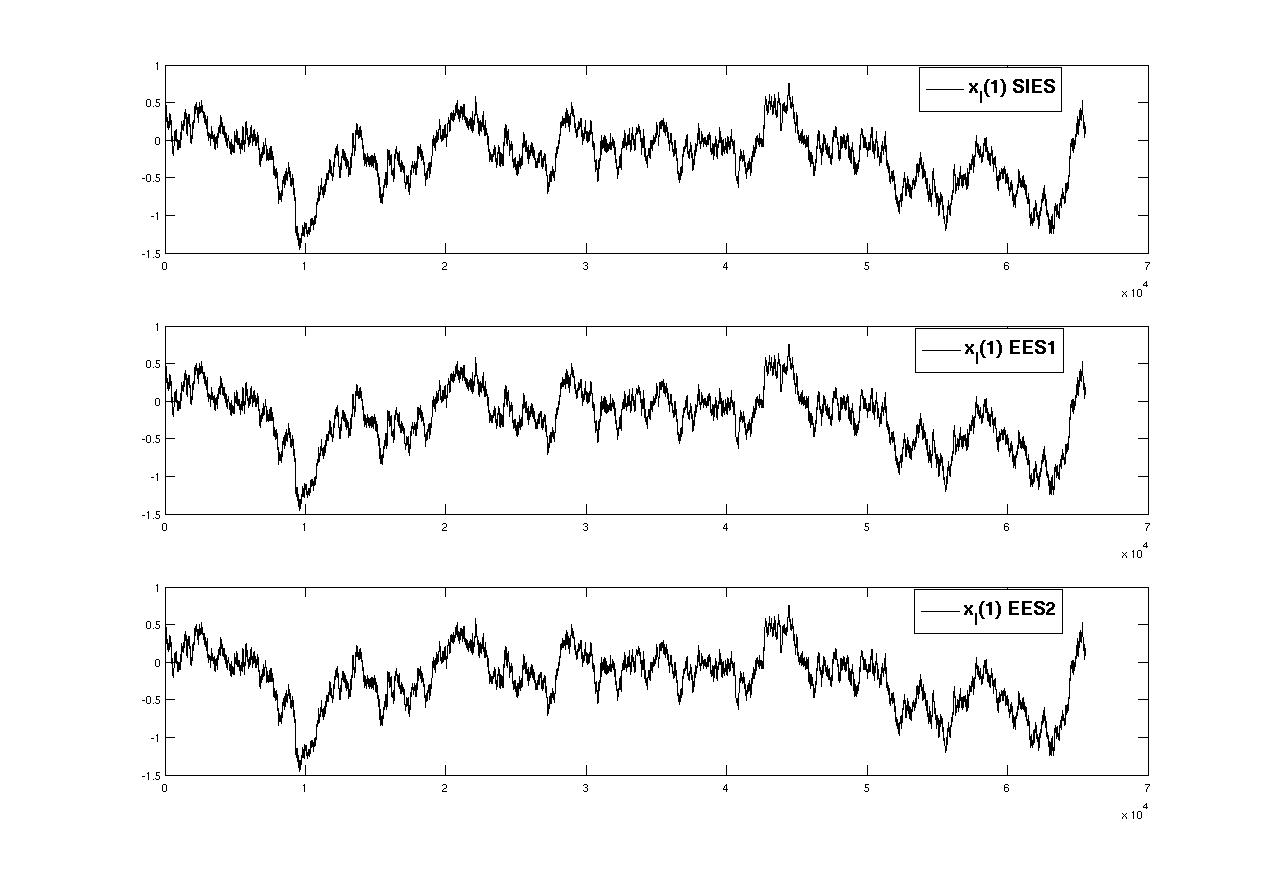}
 \caption{The first component of the chains SIES, EES1 and EES2  for $l=16$.}

\end{figure}
\newpage
%%%%%%%%%%%%%%%%%%%%%%%%%%%%%%%%%%%%%%%

\subsection{The Cost of each sheme using MC}  
We approximate for each scheme $\x(T)$ by $\x_L(sc,2^L)$ with $L=16$, and 
we look for the optimal level $lopt$ and the optimal sample size $Nopt$ such that 
\be 
MSE:=\Eb\Big[\|\Eb[\x_L(sc,2^L)]-\frac{1}{Nopt}\sum_{k=1}^{Nopt}\x_l^k(sc,2^l)\|^2\Big]=\eta^2.
\ee 
We need for $l=l_s, \ldots, L-2$ to calculate $e(sc,\Delta t_l):=\|\Eb[\x_L(sc,2^L)]-\Eb[\x_l(sc,2^l)]\|^2$. Using Monte-Carlo with $N=1000$, we obtain 
by 
\be 
e(sc,\Delta t_l)\approx \|\frac{1}{N}\sum_{k=1}^{N}\x_L^k(sc,2^L)-\frac{1}{N}\sum_{k=1}^{N}\x_l^k(sc,2^l)\|^2.
\ee 
Table 1 shows the numerical values of $e(sc,\Delta t_l)$
for each scheme and for $l=5, \ldots, 13$. 
%%%%%%%%%%%%%%%%%%%%%%%%%%%%%%%%%%%%%%%%%%%%%%%%%%%%%
%\begin{center}
%\textbf{Table of MC errors}
%\end{center}
\renewcommand{\arraystretch}{1.3} %donne la distance entre les lignes%
\setlength{\tabcolsep}{0.09cm} %donne la distance entre les collones%
\begin{table}[h!]
\begin{center}
  \begin{tabular}{|c|c|c|c|c|c|c|c|c|c|c|}
   \hline
    \scriptsize{$l$} & \scriptsize{5} & \scriptsize{6} & \scriptsize{7} & \scriptsize{8} & \scriptsize{9} & \scriptsize{10} & \scriptsize{11} & \scriptsize{12} & \scriptsize{13}  \\
    \hline
 \scriptsize{$e(SIES,\Delta_{t_l})$} & \scriptsize{0.0050}& \scriptsize{0.0080} & \scriptsize{0.0071} & \scriptsize{0.0022} & \scriptsize{0.0054} & \scriptsize{0.0066} & \scriptsize{0.0056} & \scriptsize{0.0043} & \scriptsize{0.0022}  \\ 
     \hline
\scriptsize{$e(EES1,\Delta_{t_l})$} & \scriptsize{0.0380} & \scriptsize{0.0025} & \scriptsize{0.0069} & \scriptsize{0.0043} & \scriptsize{0.0016} & \scriptsize{0.0039} & \scriptsize{0.0027} & \scriptsize{0.0032} & \scriptsize{0.0022}  \\
    \hline     
\scriptsize{$e(EES2,\Delta_{t_l})$} & \scriptsize{0.0107}& \scriptsize{0.0041} & \scriptsize{0.0044} & \scriptsize{0.0042} & \scriptsize{0.0054} & \scriptsize{0.0111} &\scriptsize{0.0041} &  \scriptsize{0.0048}& \scriptsize{0.0065} 
\\
\hline 
 \end{tabular}
  \caption{Numerical values of $e(sc,\Delta t_l)$
for each scheme and for $l=5, \ldots, 13$.}
  \end{center}
\end{table}

By fixing $\eta^2\geq \max(e(sc,\Delta t_l), sc=SIES, EES1, EES2,l=5, \ldots, 13)$, the constraint 
$e(sc,\Delta t_l) \leq \eta^2$ holds for each level $l=5, \ldots, 13$. 
The optimal level $l{opt}$ is such that $3e(sc,\Delta t_l)-2e(sc,\Delta t_{l+1})\approx \eta^2$. 
Having $lopt$ we calculate 
\be 
Var_{lopt}(sc)&:=&\sum_{i=1}^p Var(x_{lopt,i}(sc,2^{lopt})),\\
&\approx &\sum_{i=1}^p\frac{1}{N}\sum_{k=1}^N\Big|x_{lopt,i}^k(sc,2^{lopt})-\frac{1}{N}
\sum_{k=1}^Nx_{lopt,i}^k(sc,2^{lopt})\Big|^2,  
\ee 
and we derive the optimal $Nopt(sc)=\frac{Var_{lopt}(sc)}{\eta^2-e(sc,\Delta t_{lopt})}$.

%%%%%%%%%%%%%%%%%%%%%%%%%%%%%%%%%%%%%%% figure MC Nopt
The Figure 5 shows how to find graphically the optimal level $lopt$. 

\begin{figure}[!ht]
  \centering

     \includegraphics[width=14cm]{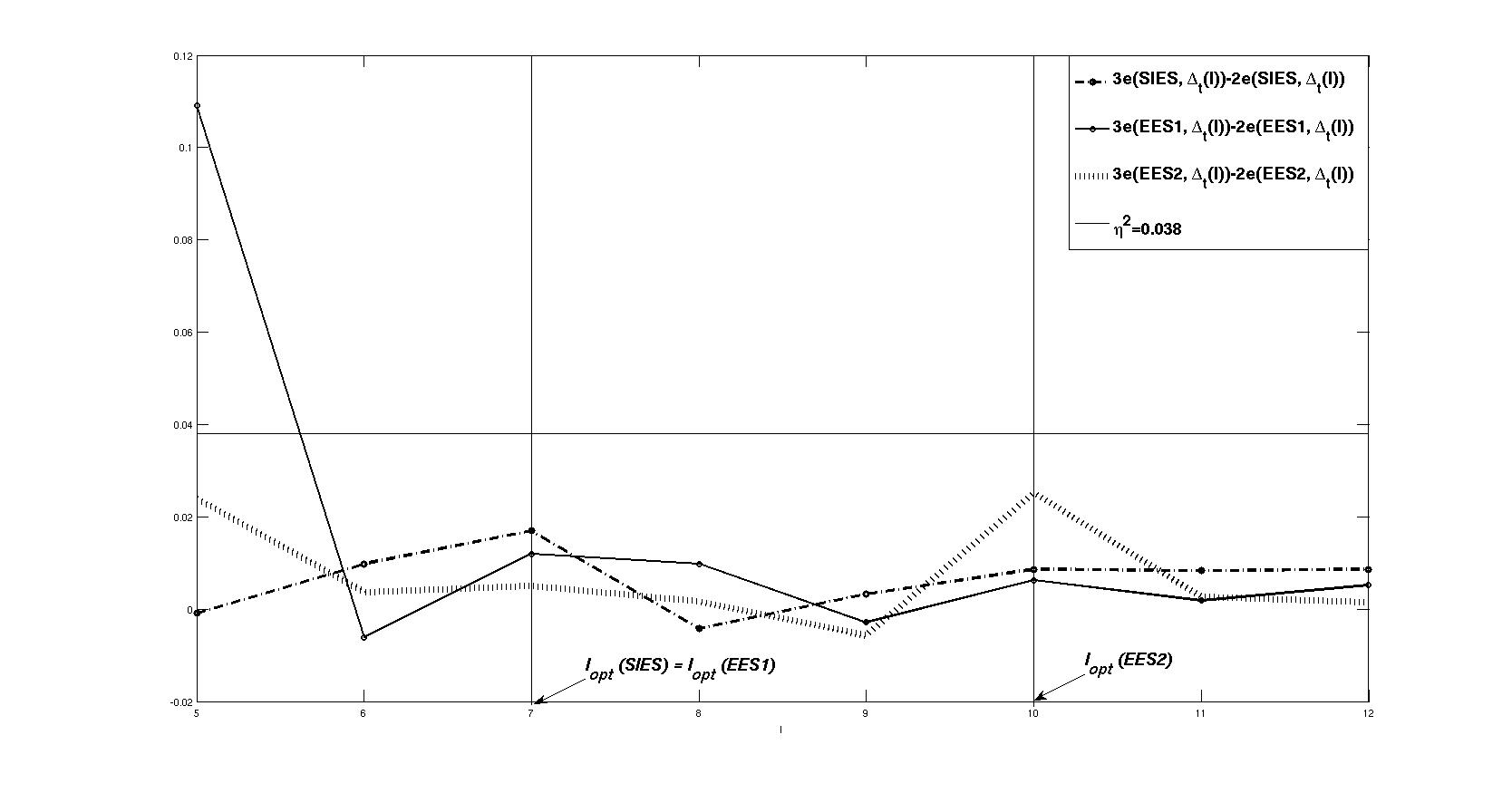}
 \caption{Graphical determination of $l_{opt}$ for the schemes SIES, EES1 and EES2.}

\end{figure}
%%%%%%%%%%%%%%%%%%%%%%%%%%%%%%%%%%%%%%%%%%%%%%

We summarize for the three schemes in the Table 2 the values of $lopt$, $Nopt$ and their cost. The scheme SIES has the lowest cost. 

\renewcommand{\arraystretch}{1} %donne la distance entre les lignes%
\setlength{\tabcolsep}{0.09cm} %donne la distance entre les collones%
\begin{table}[h!]
\begin{center}
  \begin{tabular}{|c|c|c|c|}
   \hline
   & \scriptsize{$lopt$}  & \scriptsize{$Nopt$} & \scriptsize{Cost}    \\
    \hline
 \scriptsize{SIES} &   \scriptsize{7} &   \scriptsize{70} & \scriptsize{8938}  \\ 
     \hline
\scriptsize{EES1}& \scriptsize{7} & \scriptsize{81} & \scriptsize{10427}  \\
    \hline     
\scriptsize{EES2}& \scriptsize{10} & \scriptsize{83} & \scriptsize{85035} \\
\hline 
      
\end{tabular}
 \caption{Optima level and cost of MC for each scheme.}
  \end{center}
\end{table}

%%%%%%%%%%%%%%%%%%%%%%%%%%%%%%%
\subsection{Computational cost of MLMC}  
In the Figure (6) for each scheme  
we plot $l\to e(sc,\Delta t_l)$ "equation~(\ref{el})". 
We derive graphically 
the optimal level $l_{opt}$. 

\begin{figure}[!ht]
  \centering

     \includegraphics[width=15cm]{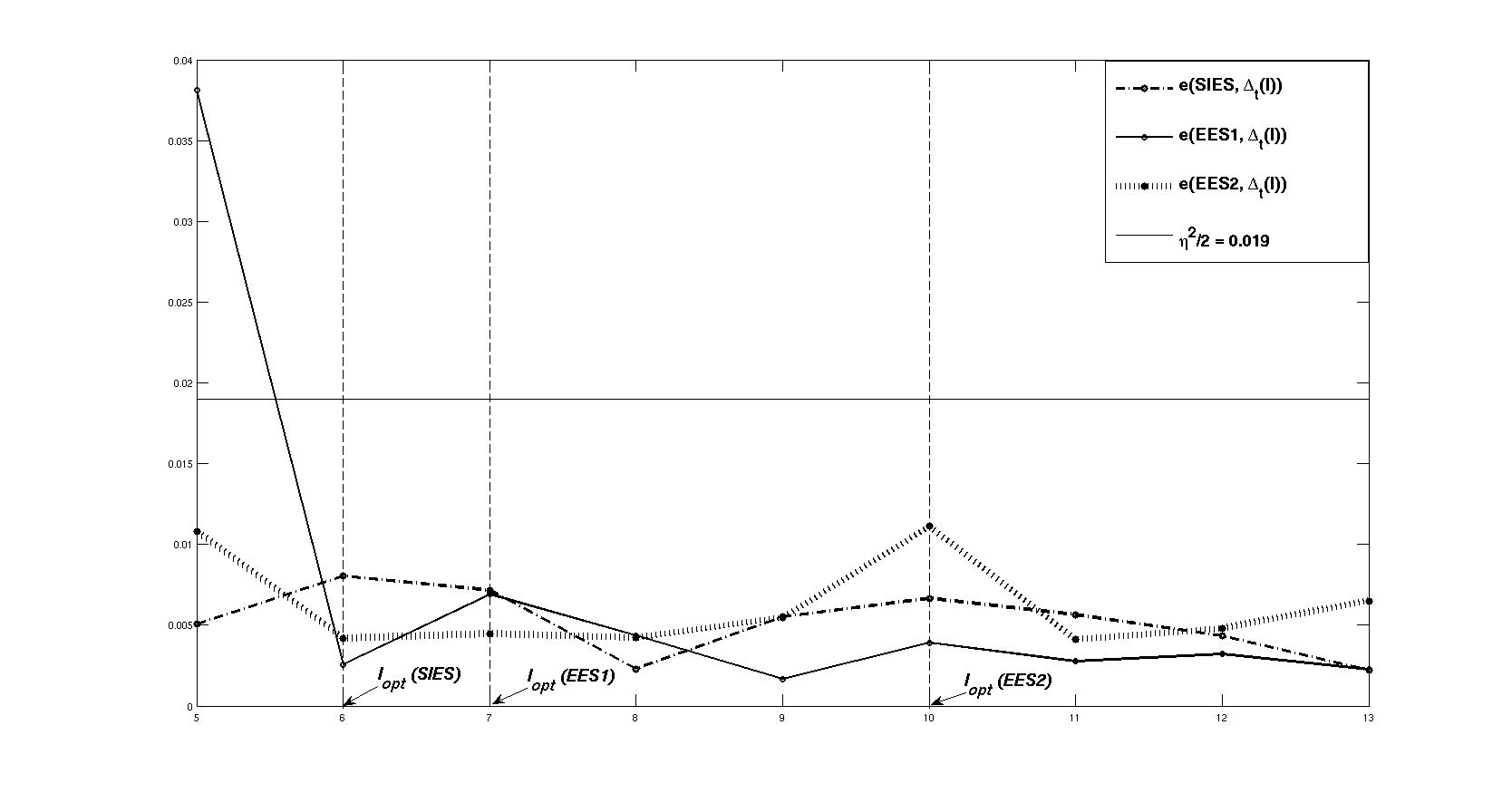}
 \caption{Graphical identification of $l_{opt}$ for each sheme.}

\end{figure}
%%%%%%%%%%%%%%%%%%%%%%%%%%%%%%%%%%%%%
We summarize for the three schemes in the Table 2 the values of $l_{opt}$, $N_{l_s}(opt), \ldots, N_{lopt}(opt)$ and their cost. 
Like MC method the scheme SIES has the lowest cost.  
\renewcommand{\arraystretch}{1} %donne la distance entre les lignes%
\setlength{\tabcolsep}{0.07cm} %donne la distance entre les collones%
\begin{table}[h!]
\begin{center}
  \begin{tabular}{|c|c|c|c|}
   \hline
      & \scriptsize{$lopt$}  & \scriptsize{$N_{5-lopt}(opt)$} & \scriptsize{Cost}    \\
    \hline
 \scriptsize{SIES} &   \scriptsize{6} &  \scriptsize{74}\quad \scriptsize{20} & \scriptsize{3639.18}  \\ 
     \hline
\scriptsize{EES1}& \scriptsize{7} &\scriptsize{132}\quad\scriptsize{40}\quad \scriptsize{16} & \scriptsize{8962.85}  \\
    \hline     
\scriptsize{EES2}& \scriptsize{10} &\scriptsize{167}\quad\scriptsize{59}\quad\scriptsize{23}\quad\scriptsize{9}\quad\scriptsize{4}\,\, \scriptsize{2} & \scriptsize{18029.47} \\
\hline 
\end{tabular}
  \caption{Optimal level and cost of MLMC for each scheme.}
  \end{center}
\end{table}
\\
\textbf{N.B.} For each $lopt$, the optimal sample sizes are 
$N_{5-lopt}:=N_{5}(opt), \ldots, N_{lopt}(opt)$,  e.g. for the scheme SIES $lopt=6$ and 
$N_{5-6}=74,20$.

%%%%%%%%%%%%%%%%%%%%%%%%%%%%%%%%%%%%%%%%%%%%%%%%%%%%%%%%%%%
\section{Markov chain Monte Carlo method MCMC} 
Using the ergodicity we suppose that the PDF of $\x(T)$ is approximated by 
$\rho(\x)=Z^{-1}\exp\big(-2\|\x\|_1-\|\A\x-\y\|^2\big)$. For the error $\eta^2$ fixed 
the cost of MCMC is the sample size  
$N$ such that 
\be 
\Eb\left[\|\Eb[\x(T)]-\frac{1}{N}\sum_{k=1}^NMCMC(k)\|^2\right]\approx \eta^2.
\ee 
Here $MCMC$ is a trajectory of the Markov Chain Monte Carlo 
having the target $\rho$. 

We recall how MCMC works. Let $k\to MC(k)$ be a Markov chain having the transition 
probability density $\pi(\x_2\,|\,\x_1) >0$ for all $\x_1, \x_2\in\Rb^p$. 
We construct from $MC$ a new  Markov chain $k\to MCMC(k)$ having the transition probability
\be 
\alpha\pi(\x_2\,|\,\x_1)d\x_2+(1-\alpha)\delta_{\x_1}(\x_2) 
\ee 
where 
\be 
\alpha=\min\Big(1,\frac{\rho(\x_2)\pi(\x_1\,|\,\x_2)}{\rho(\x_1)\pi(\x_2\,|\,\x_1)}\Big).  
\ee 
The new Markov chain $MCMC$ is ergodic and has $\rho(\x)$ as its invariant probability density function. 
We propose the Markov chains $MC(k):=\x_{lopt}(sc,k)$
for $sc=EES1, EES2$ and $MC(k)=RW(k,\sigma^2)$.
Here $RW(k,\sigma^2)$ denotes the Gaussian random walk, each step has the variance $\sigma^2$. We obtain three MCMC chains: $MCMC_{prox}(EES1)$, $MCMC_{prox}(EES2)$,
$MCMC_{RW}$. Observe that $MCMC_{prox}(EES1)$ is known as PMALA \cite{Pereyra}. 
Table 4 shows the cost of each method. 

%%%%%%%%%%%%%%%%%%%%%%%%%%%%%%%%%%%%%%%%%%%%%

\subsection{Computational cost of MCMC}  

In the table 4, we indicate the different costs of MC, $MCMC_{prox}$ and $MCMC_{RW}$. We create for each $N$, $M$ MCMC chains 
$(MCMC^i(k): k=1, \ldots, N, i=1, \ldots, M)$. 
We approximate $\Eb\left[\|\Eb[\x(T)]-\frac{1}{N}\sum_{k=1}^NMCMC(k)\|^2\right]$ by 
$\frac{1}{M}\sum_{i=1}^{M} 
\|\Eb[\x_L(sc,2^L)]-\frac{1}{N}\sum_{k=1}^NMCMC^i(k)\|^2$. 

\renewcommand{\arraystretch}{1} %donne la distance entre les lignes%
\setlength{\tabcolsep}{0.09cm} %donne la distance entre les collones%
\begin{table}[h!]
\begin{center}
  \begin{tabular}{|c|c|c|c|c|c|}
   \hline
   & \scriptsize{Cost (MC)}  & \scriptsize{Cost $(MCMC_{prox})$} & \scriptsize{Cost $(MCMC_{RW})$}   & \scriptsize{Cost $(MCMC_{RW})$}  \\
    \hline

\scriptsize{EES1}& \scriptsize{10427} & \scriptsize{5340}   & \scriptsize{3990 ( $\sigma^2=0.3$)}  & \scriptsize{17230 ( $\sigma^2=0.8$)}   \\
    \hline     
\scriptsize{EES2}& \scriptsize{85035} & \scriptsize{6200}& \scriptsize{3890 ( $\sigma^2=0.3$)}& \scriptsize{16230 ( $\sigma^2=0.8$)}   \\
\hline 
      
\end{tabular}
 \caption{Cost of MC, $MCMC_{prox}$ and $MCMC_{RW}$ for EES1 and EES2 schemes.}
  \end{center}
\end{table}

Table 4 shows that the $MCMC_{RW}$ corresponding to the proposal distribution $\Nc(0,0.3)$ is the winer. But it loses against
MLMC with the scheme SIES (see Table 2).  

{\bf Concluding remark.} In this work we 
studied the approximation of Bayesian Lasso using 
MC, MLMC and MCMC methods and three schemes Semi-implicit Euler scheme (SIES),
and two  
Explicit Euler schemes EES1 and EES2. Furthermore, we proposed a method for calculating the cost of each method and each scheme. We showed that 
the winner is MLMC with the scheme (SIES).

\end{document}